\documentclass[10pt, a4paper, english]{article}
\usepackage[margin=1.2in]{geometry}
\usepackage[english]{babel} 
\usepackage{graphicx,float,subcaption,caption} 
\usepackage{amsmath, amsthm, mathtools, booktabs, tabularx,colortbl, amsfonts,amssymb} 
\usepackage{parskip} 
\usepackage{xcolor} 
\usepackage{soul} 
\usepackage{fancyhdr} 
\usepackage{hyperref} 
\usepackage{enumitem} %
\usepackage{multicol} %
\usepackage{verbatim}
\usepackage{authblk}
\usepackage{mathtools} 
\usepackage[ruled,vlined]{algorithm2e} 
\LinesNumbered

\usepackage{bbm}
\hypersetup{
	colorlinks=true, 
	linkcolor=teal,          	
	urlcolor=blue!80!black,
    citecolor=gray,
    bookmarksdepth=4,
    pdftitle={Average distance and maximum induced forest},pdfauthor={Bouchard, Hertz}}  
\numberwithin{equation}{section} 
\numberwithin{figure}{section}
\numberwithin{table}{section} 

\newtheorem{corollary}{Corollary}[section]

\newtheorem{theorem}{Theorem}[section]
\newtheorem{property}{Property}[section]
\theoremstyle{definition}

\theoremstyle{lemma}
\newtheorem{lemma}{Lemma}[section]
\theoremstyle{conjecture}
\newtheorem{conjecture}{Conjecture}[section]\theoremstyle{formulation}

\newcommand{\diam}{\operatorname{diam}}          
\newcommand{\abs}[1]{\lvert #1\rvert}

\title{Average distance and maximum induced linear forest.}
\small\author{
	Mathieu Bouchard\textsuperscript{1} and
	Alain Hertz\textsuperscript{2}
	\\
	\footnotesize \textsuperscript{1} Perroquet AI, Québec, Canada\\
	\footnotesize mbouchard@perroquet.ai\\\vspace{4pt}
	\footnotesize \textsuperscript{2} Department of Mathematics and Industrial
	Engineering\\
	\footnotesize Polytechnique Montr\'eal - Gerad, Montr\'eal, Canada\\
	\footnotesize Corresponding author: alain.hertz@gerad.ca
}
\date{}

\allowdisplaybreaks

\begin{document}
	
	\maketitle
	\begin{abstract}
		Let $G$ be a connected graph of order $n$, let $\mu(G)$ denote its average distance, and let $L(G)$ denote the maximum order of an induced linear forest of $G$, that is, of an induced subgraph whose components are chains. Hansen, Hertz, Kilani, Marcotte and Schindl conjectured that $\mu(G) \le L(G)/2$, which can be written as $4W(G) \le L(G)\,n(n-1)$, where $W(G)$ is the Wiener index of $G$. This inequality strengthens their bound $\mu(G) \le F(G)/2$, where $F(G)$ is the maximum order of an induced forest, as well as Chung’s classical bound $\mu(G) \le \alpha(G)$. In this paper we prove the conjecture for all connected graphs.
	\end{abstract}
\setlength{\abovedisplayskip}{3pt}     
\setlength{\belowdisplayskip}{3pt}
\setlength{\abovedisplayshortskip}{3pt}
\setlength{\belowdisplayshortskip}{3pt}
\setlength{\jot}{2pt}
\section{Introduction}

The \emph{average distance} of a connected graph $G$ of order $n$, with vertex set $V_G$, denoted $\mu(G)$, is
the arithmetic mean of the distances $d_G(u,v)$ over all unordered pairs of
vertices $\{u,v\} \subseteq V_G$. For a graph $G$ of order one, we adopt the convention $\mu(G)=0$, whereas, for $|V(G)|\geq 2$,
\[
\mu(G) \;=\; \frac{1}{\binom{n}{2}} \sum_{\{u,v\}\subseteq V_G} d(u,v).
\]
This parameter has attracted considerable attention since the 1970s, notably because of its applications to communication networks and the modelling of social and biological networks, as well as its strong connections with other classical graph invariants, including the Wiener index
\cite{wiener1947}.

A natural line of research consists in bounding $\mu(G)$ from above in
terms of well-understood structural parameters. An important result in this direction was established by Chung~\cite{Chung_independence}, who proved that
the average distance of any connected graph is at most its \emph{independence
	number}~$\alpha(G)$, i.e.\ the maximum size of an independent set:
\[
\mu(G) \;\leq\; \alpha(G).
\]
This result was later strengthened by Dankelmann \cite{Dankelmann1994}, who gave a sharp upper bound on the average distance of a connected graph of given order and independence number, together with a characterization of the extremal graphs. This inequality  elegantly relates a metric invariant to a
purely combinatorial one. 

Let $B(G)$ denote the maximum number of vertices in an \emph{induced bipartite subgraph} of a graph $G$. Clearly, $B(G)\leq 2\alpha(G)$.  Around 1992, Fajtlowicz formulated the following conjecture with the help of the Graffiti system:
\[
\mu(G) \;\leq\; \tfrac{1}{2}\,B(G).
\]
It was included, along with many others, in the electronic file Written on the Wall \cite{graffiti}. More recently, Hansen et~al.~\cite{HansenEtAl2009} established a better bound: the average distance is at most half the order
of a largest \emph{induced forest} of~$G$.  Denoting by $F(G)$ the
maximum number of vertices in an induced forest of~$G$, their result
reads
\[
\mu(G) \;\leq\; \tfrac{1}{2}\,F(G).
\]
In the same paper, the authors raised the question of whether the bound
could be strengthened by further restricting the class of induced
subgraphs considered.  Specifically, they conjectured that
an analogous inequality holds for \emph{induced linear forests}, that
is, induced subgraphs in which every connected component is a chain.
Letting $L(G)$ denote the maximum order of an induced linear
forest of~$G$, the conjecture states:

\begin{conjecture}\label{conj:main}
	For every connected graph $G$,
	\[
	\mu(G) \;\leq\; \tfrac{1}{2}\,L(G).
	\]
\end{conjecture}

Since every induced linear forest is in particular an induced forest,
we have $L(G) \leq F(G)$ for every graph~$G$, so
Conjecture~\ref{conj:main} is strictly stronger than the result of
Hansen et~al.  The present paper confirms this conjecture in full
generality. The best previous result, due to DeLaViña and Waller \cite{DeLaVinaWallerNote}, was $\mu(G) < \frac{L(G)+1}{2}$. 

\begin{figure}[H]
	\centering
	\includegraphics[width=0.1\linewidth]{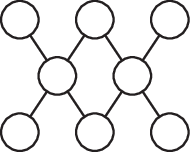}
	\caption{A graph $G$ for which $\mu(G)$, $\alpha(G)$, $\frac{1}{2}B(G)$, $\frac{1}{2}F(G)$ and $\frac{1}{2}L(G)$ all take distinct values.}
	\label{fig:differences}
\end{figure}

The graph depicted in Figure \ref{fig:differences}, which was previously considered in \cite{HansenEtAl2009}, provides an illustrative example in which all the upper bounds under consideration yield different values. Moreover, none of these bounds is attained, showing that each of them can be strict, even though they provide distinct estimates for the same graph. Indeed, for this graph, we have $$\mu(G)=\frac{15}{7}<3=\frac{1}{2}L(G)<\frac{7}{2}=\frac{1}{2}F(G)<4=\frac{1}{2}B(G)<6=\alpha(G).$$
For two vertices $u$ and $v$ in a graph $G$, let $d_G(u,v)$ denote the distance between $u$ and $v$, and let 
$$
W(G)=\sum_{\{u,v\}\subseteq V_G}d_G(u,v).
$$
Hence, $\mu(G)=\frac{2W(G)}{n(n-1)}$ and Conjecture \ref{conj:main} is therefore equivalent to $L(G)n(n-1)\geq 4W(G)$ for every connected graph $G$ of order $n$.

A \emph{dominating set} in a graph $G$ is a set $D\subseteq V_G$ of vertices such that every vertex in $V_G\setminus D$ is adjacent to at least one vertex in $D$. A connected dominating set of $G$ is a dominating set $D$ such that $G[D]$ is connected. We denote by $\gamma_c(G)$ the minimum number of vertices in a connected dominating set of $G$.

The \emph{diameter} of a connected graph $G$, denoted by $\diam(G)$, is the maximum distance between any two vertices of $G$, that is,
$$\diam(G)=\max_{\{u,v\}\subseteq V_G}d_G(u,v).$$
We prove Conjecture \ref{conj:main} by considering several cases according to the values of $\diam(G)$ and $\gamma_c(G)$. More precisely, the paper is organized as follows. Section \ref{sec:2} introduces some definitions and properties. In Section \ref{sec:3}, we prove Conjecture \ref{conj:main},by induction on $n$, in two simple cases, namely when $L(G)\leq 3$ and when $L(G)\leq \diam(G)+2$. In Section \ref{sec:4}, we establish a relationship between $L(G)$ and $\gamma_c(G)$ and, in particular, show that $L(G)\geq \gamma_c(G)+1$. Section \ref{sec:5} establishes the conjecture for $L(G)=\gamma_c(G)+1$. Finally, Section \ref{sec:6} addresses the remaining case, where $L(G)\geq \gamma_c(G)+2$ and $L(G)\geq \diam(G)+3$.

\section{Notations and properties}\label{sec:2}
For a graph $G$, we denote its vertex set by $V_G$ and its edge set by $E_G$. For a subset $W\subseteq V_G$ of vertices of $G$, we denote by $G[W]$ the subgraph of $G$ induced by $W$, and by $N^1_G(W)$ the set of vertices in $V_G\setminus W$ that are adjacent to at least one vertex in $W$. Also, for $W\subseteq V_G$, we denote by $N^2_G(W)$ the set of vertices of $G$ whose minimum distance to a vertex of $W$ is equal to 2, meaning that $N^2_G(W)=N^1_G(W\cup N^1_G(W))=N^1_G(N^1_G(W))\setminus W$.

For a connected graph $G$ of order $n$, we define 
$$
\delta(G)=L(G)n(n-1)-4W(G).$$
Thus, proving Conjecture \ref{conj:main} amounts to showing that $\delta(G)\geq 0$ for every connected graph $G$ of order $n$.
For a vertex $u$ in a connected graph $G$, we denote $$w_G(u)=\sum_{v\in V_G}d_G(u,v)\quad\quad\mbox{and}\quad\quad\tau_G(u)=L(G)(n-1)-2w_G(u).
$$
Since $\sum_{u\in V_G}w_G(u)=2W(G)$, it follows that
$$\sum_{u\in V_G}\tau_G(u)=\delta(G).$$
Moreover, given any two distinct vertices $u$ and $v$, we have
\begin{equation}\label{eq:tau-sum}
	\tau_G(u)+\tau_G(v)=2\Big(L(G)(n-1)-\sum_{w\in V_G}\big(d_G(u,w)+d_G(v,w)\big)\Big).
\end{equation}

For a tree $T$ and an edge $e\in E_T$, let $a(T,e)$ denote the asymmetry of $T-e$, defined as the difference between the orders of the largest and smallest connected components of $T-e$. Also, let $\pi(T,e)$ denote the product of the orders of the two connected components of $T-e$. The following property is a well-known theorem due to Wiener \cite{wiener1947} .  
\begin{property}  \label{pro:Wiener}
	If $T$ is a tree, then $W(T)=\sum_{e\in E_T}\pi(T,e)$.
\end{property}
We define, for a graph $H$ and an integer $h$, $$\underline{W}(H,h)=\sum_{\{u,v\}\subseteq V_H}
\min \{d_H(u,v),h\}.$$ 
As we show below, if $H$ is a spanning subgraph of $G$, and if $h$ is at least $\diam(G)$, then $\underline{W}(H,h)$ is an upper bound on $W(G)$, while it is clearly a lower bound on $W(H)$.
\begin{property}\label{prop:W}
Let $H$ be a connected spanning subgraph of $G$ (not necessarily induced), where $G$ is connected, and let $h\geq \diam(G)$. Then $W(G)\leq \underline{W}(H,h)\leq W(H)$.
\end{property}
\begin{proof}
	By definition, $\underline{W}(H,h)\leq W(H)$. Moreover, for every pair of vertices $u,v$ of  $G$, we have 
	\begin{itemize}
		\item $d_G(u,v)\leq d_H(u,v)$ because $H$ is a spanning subgraph of $G$;
		\item $d_G(u,v)\leq \diam(G)\leq h$. 
	\end{itemize}
	Hence, for every pair of vertices $u,v$, we have $d_G(u,v)\leq \min\{d_H(u,v),h\}$. Summing over all pairs of vertices $u,v$ yields $W(G)\leq \underline{W}(H,h)$.
\end{proof}

\begin{property}\label{pro:delta}
If $u$ is a vertex in a connected graph $G$ such that $G-u$ is also connected, then $\delta(G)-\delta(G-u)\geq 2\tau_G(u)$.
\end{property}
\begin{proof}Every induced linear forest of $G-u$ is an induced linear forest of $G$, so $L(G-u)\le L(G)$. Distances in $G-u$ are at least those in $G$, so $W(G-u)\ge W(G)-w_G(u)$. Therefore
\begin{align*}
	\delta(G)-\delta(G-u)
	&=L(G)n(n-1)-4W(G)-L(G-u)(n-1)(n-2)+4W(G-u)\\
	&\ge L(G)n(n-1)-4W(G)-L(G)(n-1)(n-2)+4W(G)-4w_G(u)\\
	&=2L(G)(n-1)-4w_G(u)=2\tau_G(u).\qedhere
\end{align*}
\end{proof}
\begin{corollary}\label{cor:delta}
	Let $G$ be a connected graph of order $n$, and suppose that Conjecture \ref{conj:main} holds for all connected graphs of order at most $n-1$. If $G$ contains a vertex $u$ such that $G-u$ is connected and $\tau_G(u)\geq 0$, then the conjecture also holds for  $G$.
	\end{corollary}
	\begin{proof}
It follows from Property \ref{pro:delta} that $
	\delta(G){-}\delta(G-u)
	{\geq} 2\tau_G(u)$.
Since $\tau_G(u)\geq 0$, we have $\delta(G)\geq \delta(G-u)$. As $G-u$ is connected, the induction hypothesis yields $\delta(G-u)\geq 0$. Hence, $\delta(G)\geq 0$.
\end{proof}

\begin{figure}[H]
	\centering
	\includegraphics[width=0.63\linewidth]{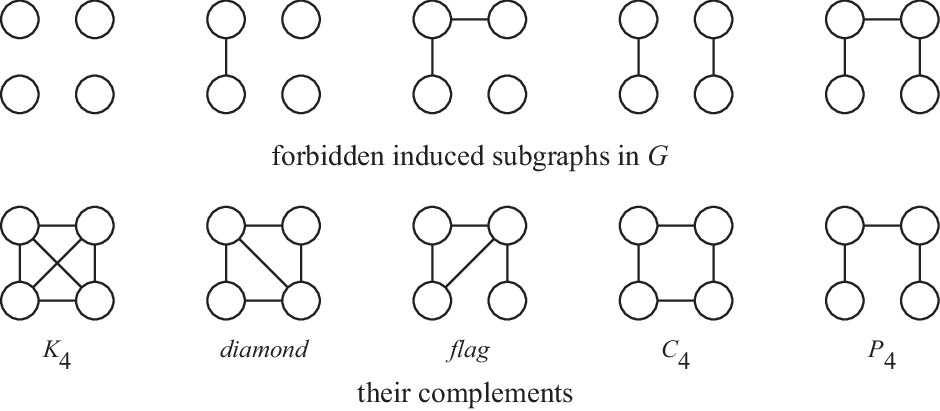}
	\caption{Forbidden induced subgraphs of graphs $G$ with $L(G)=3$ and their  complements.}
	\label{fig:C4P4}
\end{figure}

The following property provides an upper bound on the number of missing edges in a graph $G$ satisfying $L(G)\leq 3$.
\begin{property}\label{prop:nonedge}
	Every graph $G$ satisfying $L(G)\leq 3$ has at most $|V_G|$ nonadjacent vertex pairs.
\end{property}
\begin{proof}
	If $L(G)\leq 2$, then either $G$ is a clique, with no missing edges, or it consists of exactly two isolated vertices, with one missing edge. Thus, in both cases, the number of missing edges is strictly less than $|V_G|$.
		So, assume $L(G)=3$. The five graphs at the top of Figure~\ref{fig:C4P4} are the only minimal forbidden induced subgraphs of $G$. Their complements, shown at the bottom of the same figure, are therefore forbidden induced subgraphs of the complement $\overline{G}$ of $G$. Since $P_4$ and $C_4$ are forbidden in $\overline{G}$, it is well known \cite{Golumbic1978} that every connected component of $\overline{G}$ contains a universal vertex, that is, a vertex adjacent to every other vertex in the component.
	
	Consider a connected component of $\overline{G}$ and let $u$ be a universal vertex in this component. If the component contains two adjacent vertices $v$ and $w$, both distinct from $u$, then it cannot contain a fourth vertex, as this would induce a $K_4$, a diamond, or a flag in $\overline{G}$ (see Figure \ref{fig:C4P4}). Hence, every connected component of $\overline{G}$ is either a triangle or a star. It follows that $\overline{G}$ has at most $|V_G|$ edges, with equality attained when $\overline{G}$ is a disjoint union of triangles. We conclude that $G$ has at most $|V_G|$ nonadjacent vertex pairs. 
\end{proof}
\section{The cases $L(G)\leq 3$ and $L(G)-\diam(G)\leq 2$}\label{sec:3}
We begin this section by showing that Conjecture \ref{conj:main} holds for connected graphs $G$ with $L(G)\leq 3$.

\begin{theorem}\label{thm:L3}
	Every connected graph $G$ of order $n$ with $L(G){\leq} 3$ satisfies $L(G)n(n-1){\geq} 4W(G)$.
\end{theorem}

\begin{proof}
	If $G$ consists of a single vertex, then $L(G)n(n-1)=4W(G)=0$, and if $G$ is complete of order $n\geq 2$, then  $L(G)n(n-1)=4W(G)=2n(n-1)$. Hence, we may assume that $G$ is a non-complete connected graph. As any induced chain in $G$ joining two non-adjacent vertices is a linear forest, we have $L(G)\geq 3$. Combined with the assumption that $L(G)\leq 3$, this yields $L(G)=3$. Consequently, $d_G(u,v)\leq 2$ for every pair of vertices $u,v\in V_G$, implying that every pair of non-adjacent vertices is at distance $2$, and
$$
	W(G)=|E_G|+2\Big(\binom n2-|E_G|\Big)=n(n-1)-|E_G|.
$$
	If $G$ has order $n=3$, then its connectedness and non-completeness imply that it has exactly two edges, yielding $L(G)n(n-1)=18>16=4W(G)$.
		We may therefore assume that $n\geq 4$. 
		
		Every induced subgraph $H$ of $G$ on four vertices contains at least three edges, since otherwise $L(G)\geq L(H)\geq 4$. Let $\mathcal{S}$ be the family of all four-vertex subsets of $V_G$. Each edge of $G$ belongs to exactly $\binom{n-2}2$ sets $S$ in $\mathcal{S}$. Hence, 
	$$
	|E_G|\binom{n-2}2=\sum_{S\in \mathcal{S}}|E_{G[S]}|\ \ge\ 3\binom n4.$$
	
	This implies that $|E_G|\geq \frac{3\binom n4}{\binom{n-2}2}=\frac{n(n-1)}4$, which yields $W(G)\leq \frac{3}{4}n(n-1)$. Consequently, $L(G)n(n-1)=3n(n-1)\geq 4W(G)$.
\end{proof}

Since every chain joining two vertices $u$ and $v$ at distance $\diam(G)$ is a linear forest, we have $L(G)\geq \diam(G)+1$ for every connected graph $G$. We now prove that the Conjecture \ref{conj:main} holds when $L(G)=\diam(G)+1$, assuming that it holds for all connected graphs of smaller order.

\begin{theorem}\label{thm:diam1}
	Let $G$ be a connected graph of order $n$ satisfying $L(G)=\diam(G)+1$. If Conjecture \ref{conj:main} holds for all connected graphs of order at most $n-1$, then it also holds for  $G$.
\end{theorem}
\begin{proof}
	As shown in Theorem \ref{thm:L3}, the conjecture holds when $L(G)\leq 3$. We may therefore assume that $L(G)\geq 4$, which implies that $\diam(G)\geq 3$. Let $u_1$ and $u_2$ be vertices at distance $\diam(G)$ in $G$, and let $P$ be a shortest chain connecting them. Hence, $|V_P|=L(G)$, and thus $V_G=V_P\cup N^1_G(P)$. 
	
	Let $Q$ be the set of vertices in $N^1_G(P)$ that have exactly one neighbor on $P$. Note that no vertex of $Q$ is adjacent to either endpoint of $P$, since adding such a vertex to $P$ would yield a chain (and hence a linear forest) on $L(G)+1$ vertices. For any two vertices $q_1$ and $q_2$ in $Q$, with respective neighbors $p_1$ and $p_2$ on $P$, we have $d_P(p_1,p_2)\geq 3$. Indeed, \begin{itemize}
		\item	if $d_P(p_1,p_2)\leq 1$, then $(V_P\cup\{q_1,q_2\})\setminus\{p_1\}$ induces a linear forest of order $L(G)+1$.
		\item	if $d_P(p_1,p_2)=2$, then  there is a vertex $p_3$ on $P$ adjacent to both $p_1$ and $p_2$, and  $(V_P\cup\{q_1,q_2\})\setminus\{p_3\}$ is a linear forest of order $L(G)+1$.
\end{itemize}
Therefore, $|Q|\leq \lfloor\frac{\diam(G)+1}{3}\rfloor$.
Consider the quantity $d_G(u_1,v)+d_G(u_2,v)-\diam(G)$ for each vertex $v$ of $G$. It is equal to $0$ when $v\in V_P$, at most $1$ when $v\in N_G^1(P)\setminus Q$, and at most $2$ when $v\in Q$. The sum of this quantity over all vertices of $G$ is therefore at most $n-\diam(G)-1+\lfloor\frac{\diam(G)+1}{3}\rfloor$.

Let $H$ be the graph obtained from $G$ by deleting $u_1$ and $u_2$.
Note that $G-u_1$ is connected, since otherwise, adding to $P$ a vertex outside $P$ that is adjacent to $u_1$ would yield a linear forest with $L(G)+1$ vertices. Thus, $H=(G-u_1)-u_2$ is also connected. Indeed, if $H$ were disconnected, then $u_2$ would have a neighbor $v$ outside $P$. Moreover, $v$ would necessarily be adjacent to $u_1$, since otherwise adding $v$ to $P$ would yield a linear forest with $L(G)+1$ vertices. Hence, we would have $\diam(G)=d_G(u_1,u_2)\leq 2$, a contradiction.

Every linear forest of $H$ is also a linear forest of $G$, so $L(H)\leq L(G)$. Moreover, given that $d_H(x,y)\geq d_G(x,y)$ for every pair of vertices $x,y\in V(H)$, we have $$W(H)\geq W(G)-w_G(u_1)-w_G(u_2)+d_G(u_1,u_2)=W(G)-w_G(u_1)-w_G(u_2)+\diam(G).$$
It follows that
\begin{align*}
	\delta(H)&=L(H)(n{-}2)(n{-}3)-4W(H)\\
	&\leq L(G)(n{-}2)(n{-}3)-4\Big(W(G)-w_G(u_1)-w_G(u_2)+\diam(G)\Big)\\
	&=L(G)(n{-}2)(n{-}3)-4W(G)+4\Big(\sum_{v\in V_G}\big(d_G(u_1,v){+}d_G(u_2,v){-}\diam(G)\big)+\diam(G)(n{-}1)\Big)\\[-4pt]
	&\leq L(G)(n{-}2)(n{-}3)-4W(G)+4\Big(n-\diam(G)-1+\lfloor\frac{\diam(G)+1}{3}\rfloor+\diam(G)(n{-}1)\Big).
\end{align*}

Using $\delta(G)=L(G)n(n-1)-4W(G)$ and $L(G)=\diam(G)+1$, we obtain 
\begin{align*}
	\delta(G)-\delta(H)&\geq (\diam(G)+1)(4n-6)-4\Big(n-\diam(G)-1+\lfloor\frac{\diam(G)+1}{3}\rfloor+\diam(G)(n-1)\Big)\\
	&=2\Big(\diam(G)-1-2\lfloor\frac{\diam(G)+1}{3}\rfloor\Big).
\end{align*}
As $\diam(G){\geq} 3$, we have $\diam(G)-1\geq 2\lfloor\frac{\diam(G)+1}{3}\rfloor$, yielding $\delta(G)\geq \delta(H)$. Furthermore, since $H$ is connected and has $n-2$ vertices, the induction hypothesis yields $\delta(H)\geq 0$, and thus $\delta(G)\geq 0$.
\end{proof}

To conclude this section, we prove the same result as in the preceding theorem, under the assumption that $L(G)=\diam(G)+2$.

\begin{theorem}\label{thm:diam2}
	Let $G$ be a connected graph of order $n$ satisfying $L(G)=\diam(G)+2$. If Conjecture \ref{conj:main} holds for all connected graphs of order at most $n-1$, then it also holds for  $G$.
\end{theorem}
\begin{proof}
	First, observe that $L(G)=\diam(G)+2$ implies that $G$ is not complete, and hence that $\diam(G)\geq 2$. Let $u$ and $v$ be vertices at distance $\diam(G)$ in $G$, and let $P$ be a shortest chain connecting them. There is at most one vertex in $G$ outside $V_P\cup N^1_G(P)$, since otherwise adding two such vertices to $V_P$ would induce a linear forest of $G$ of order $\diam(G)+3>L(G)$.
	Moreover, since $G$ is connected, the possible vertex outside $V_P\cup N^1_G(P)$ belongs to $N^2_G(P)$.
	
	Let us examine the quantity $d_G(u,w)+d_G(v,w)-\diam(G)$ for every vertex $w$ of $G$. It equals $0$ for $w\in V_P$, is at most $2$ for $w\in N_G^1(P)$, and is at most $4$ for $w\in N_G^2(P)$.
	Therefore, substituting $\diam(G)+2$ for $L(G)$ in  (\ref{eq:tau-sum}), we obtain
	\begin{align*}
		\tau_G(u)+\tau_G(v)&=2\Big((\diam(G)+2)(n-1)-\sum_{w\in V_G}\big(d_G(u,w)+d_G(v,w)-\diam(G)\big)-n\diam(G)\Big)\\[-4pt]
		&\geq 2\Big((\diam(G)+2)(n-1)-\big(2(n-\diam(G)-1)+2\big)-n\diam(G)\Big)\\&=2(\diam(G)-2)\geq 0.\end{align*}
	
	Hence, at least one of $\tau_G(u)$ and $\tau_G(v)$ is nonnegative. Moreover, removing either endpoint of $P$ from $G$ leaves a connected graph; otherwise, a neighbor of that endpoint outside $P$ would be at distance $\diam(G)+1$ from the other endpoint of $P$, contradicting the definition of the diameter. Therefore, Corollary \ref{cor:delta} implies that the conjecture also holds for $G$.
\end{proof}

\section{Linear Forests and Connected Dominating Sets}\label{sec:4}

The following algorithm heuristically constructs a linear forest $F$ in a graph $G$. It first determines a maximal chain $P_0$ and sets $F{=}P_0$. As long as $N^2_G(V_F)$ is nonempty, additional chains are iteratively added to $F$. More precisely, when $F$ contains $i$ chains, a vertex $v_{i}$ is selected from $N^2_G(V_F)$, and a maximal chain $P_{i}$ containing $v_{i}$ is determined in $G[V_G\setminus (V_F\cup N^1_G(V_F))]$ and added to $F$. Since $v_{i}\in N^2_G(V_F)$, there exists a vertex in $N^1_G(V_F)$, which we denote by $d_i$, that is adjacent to $v_{i}$. Note that although the vertices $v_i$ are all distinct, this is not necessarily the case for the vertices $d_i$.

\begin{algorithm}[h!]\label{alg:Diamond-Decycling}
	
	\SetKwBlock{Input}{Input}{}
	\SetKwBlock{Output}{Output}{}
	\SetKwBlock{Main}{Main}{}
	\Input{
		A graph $G$.
	}
	\Output{
		A linear forest $F$ of $G$.
	}
	\Main{
		Determine a maximal chain $P_0$ in $G$ and set $F\gets P_0$ and $i\gets 1$.\\
		\While{ $N^2_G(V_F)\neq \emptyset$}
		{
			Choose a vertex $v_i$ in $N^2_G(V_F)$ as well as a vertex $d_i\in N^1_G(V_F)$ adjacent to $v_i$.\\
			Choose a maximal chain $P_i$ in $G[V_G\setminus (V_F\cup N^1_G(V_F))]$ containing $v_i$.\\
			Set $F\gets F\cup P_i$ and $i\gets i+1$.
		}
		\Return{$F$}
	}
	\caption{\texttt{Build-Linear-Forest}}
\end{algorithm}
\begin{figure}[H]
	\centering
	\includegraphics[width=0.21\linewidth]{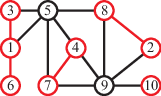}
	\caption{Illustration of the \texttt{Build-Linear-Forest} heuristic.}
	\label{fig:BLF}
\end{figure}
An illustration of the algorithm is given in Figure~\ref{fig:BLF}. Whenever several choices are possible, we select the vertex with the smallest index. The first chain is therefore constructed by starting with vertex 1 and successively adding vertices 3 and 6. At this point, $N^2_G(V_F)$ consists of vertices 4, 7, and 8. We select $v_1=4$, for which $d_1=5$, and add vertex 7 to vertex 4, thereby forming the chain $P_1$, which is then added to $F$.
The set $N^2_G(V_F)$ now consists of vertices 2, 8, and 10. We select $v_2=2$, for which $d_2=9$, and add vertex 8 to vertex 2 to form the chain $P_2$. After adding $P_2$ to $F$, only vertex 10 remains in $N^2_G(V_F)$. We thus have $v_3=10$ and $d_3=9$. The algorithm terminates with the linear forest shown in red.

Using this algorithm, we can now prove that $\gamma_c(G)\leq L(G)-1$. 
This result was previously established as Theorem 4 in an unpublished research report by DeLaVi{\~n}a and Waller \cite{DeLaVinaWallerNote}, using a construction similar to the \texttt{Build-Linear-Forest} algorithm.\\

\begin{theorem}\label{thm:bound}
	For every connected graph $G$ of order at least 2, we have
	$\gamma_c(G)\leq L(G)-1.$
\end{theorem}
\begin{proof}
	Let $P_0,P_1,\ldots,P_r$ be the chains of a linear forest $F$ generated by the \texttt{Build-Linear-Forest} algorithm, and let $n_i$ denote the number of vertices in chain $P_i$, for $i=0,1,\ldots,r$. We therefore have $
	L(G)\geq |V_F|=\sum_{i=0}^r n_i.
	$
	
	Let $p_i^1,p_i^2,\ldots,p_i^{n_i}$ denote the vertices of $P_i$, where $p_i^j$ is adjacent to $p_i^{j+1}$ for $j=1,2,\ldots,n_i-1$. Let $w_0$ be an endpoint of $P_0$, that is, $w_0=p_0^1$ or $w_0=p_0^{n_0}$. For $i=1,2,\ldots,r$, we choose a vertex $w_i$ of $P_i$ as follows. If $n_i=1$, we set $w_i=v_i$. If $n_i\geq 2$ and $v_i$ is an endpoint of $P_i$, we let $w_i$ be the other endpoint of $P_i$. Otherwise, we let $w_i$ be any endpoint of $P_i$. We now prove that the set obtained by removing $w_0,w_1,\ldots,w_r$ from $V_F$ and adding $d_1,d_2,\ldots,d_r$ is a connected dominating set of $G$.
	
	So let $ D=\left(V_F\setminus\{w_0,w_1,\ldots,w_r\}\right)\cup\{d_1,d_2,\ldots,d_r\}.$ We first prove that $D$ is a dominating set. Consider a vertex $u$ of $G$ that does not belong to $D$.
	\begin{itemize}
		\item If $u$ is one of the vertices $w_i$, there are two possibilities. If $w_i=v_i$ for some $i\geq1$, then $P_i$ consists of a single vertex, so $w_i$ is adjacent to $d_i\in D$. Otherwise, $w_i$ has a neighbor in $P_i$ that belongs to $D$.
		\item If $u$ is not one of the vertices $w_i$, then $u$ is adjacent to at least one vertex of $F$. At least one of these neighbors does not belong to $\{w_0,w_1,\ldots,w_r\}$. Indeed, suppose that all neighbors of $u$ in $F$ are among the vertices removed to obtain $D$, and let $i$ be the smallest index such that $u$ is adjacent to $w_i$. Then $P_i$ cannot be a maximal chain, since $u$ could be added to $P_i$, contradicting the construction rules of \texttt{Build-Linear-Forest}.
		\end{itemize}
		
		It remains to prove that $G[D]$ is connected. For $i=0,1,\ldots,r$, let $F_i$ denote the linear forest induced by the chains $P_0,P_1,\ldots,P_i$, and let $
		D_i=\left(V_{F_i}\setminus\{w_0,w_1,\ldots,w_i\}\right)\cup\{d_1,d_2,\ldots,d_i\}
		$.
		We prove by induction on $i$ that $G[D_i]$ is connected. Since $D_r=D$, this will establish that $G[D]$ is connected. For $i=0$, $G[D_0]$ is a chain containing at least one vertex and is therefore connected. Suppose that $G[D_{i-1}]$ is connected for some $i-1<r$. We show that $G[D_i]$ is also connected. 
		
		\begin{itemize}
			\item 
		By construction, $d_i$ is adjacent to at least one vertex of $F_{i-1}$. If all its neighbors in $F_{i-1}$ belong to $\{w_0,w_1,\ldots,w_{i-1}\}$, let $j$ be the smallest index such that $d_i$ is adjacent to $w_j$. Then $d_i$ could be added to $P_j$, contradicting the maximality of $P_j$ in the construction of \texttt{Build-Linear-Forest}. Thus, $d_i$ has a neighbor in $D_{i-1}$, and $D_{i-1}\cup\{d_i\}$ induces a connected subgraph of $G$.
		\item If $n_i=1$, then $D_i=D_{i-1}\cup\{d_i\}$, and $G[D_i]$ is therefore connected. Otherwise, $v_i\in D_i$, and $G[D_i]$ is connected because $v_i$ is adjacent to $d_i\in D_i$, while every other vertex of $P_i$ that belongs to $D_i$ is connected to $v_i$ by a chain whose vertices all belong to $D_i$.
	\end{itemize}
	This completes the proof, since $D$ is a connected dominating set containing at most $|V_F|-1$ vertices, which implies
	$
	\gamma_c(G)\leq |D| \leq |V_F|-1\leq L(G)-1.
	$
\end{proof}

\begin{theorem}\label{thm:dominating}
Let $G$ be a connected graph of order at least $2$. If $\gamma_c(G)=L(G)-1$, then every maximum induced chain of $G$ is a connected dominating set of $G$.
\end{theorem}
\begin{proof}
	Suppose that $\gamma_c(G)=L(G)-1$ and that some maximum induced chain in $G$ is not dominating. We will show that this leads to a contradiction. Since $\gamma_c(G)=L(G)-1=1$ when $G$ is a clique, in which case every maximum induced chain consists of two adjacent vertices and is dominating, $G$ cannot be a clique. Thus, every maximum induced chain in $G$ contains at least three vertices.	
	Consider an execution of \texttt{Build-Linear-Forest} in which the first chain $P_0$ constructed at Step 6 is a maximum induced chain of $G$ that is not dominating. Let $p_0^1,p_0^2,\ldots,p_0^{n_0}$ denote the vertices of $P_0$, ordered so that $p_0^j$ is adjacent to $p_0^{j+1}$ for $j=1,2,\ldots,n_0-1$. The algorithm must then construct at least one additional chain, which means that a vertex $d_1$ is selected at Step 8. Among all possible executions of \texttt{Build-Linear-Forest}, choose one for which $P_0$ is a maximum induced chain of $G$ that is not dominating, and the largest index $i$ such that $d_1$ is adjacent to $p_0^i$ is minimized. Let $F$ be the linear forest generated by this execution. Let $D$ be the connected dominating set constructed from $F$, with $w_0=p_0^{n_0}$.
 We therefore have
	$
	\{p_0^1,p_0^2,\ldots,p_0^{n_0-1},d_1\}\subseteq D
	$, with $n_0\geq 3$.
	It follows from the proof of Theorem~\ref{thm:bound} that
	$$ \gamma_c(G)\leq |D|\leq |V_F|-1\leq L(G)-1=\gamma_c(G).$$
	Thus, equality holds throughout, implying that $D$ is a minimum connected dominating set. Since $p_0^1\in D$, removing this vertex from $D$ yields a set that is no longer a connected dominating set. We now show that $G[D\setminus\{p_0^1\}]$ is connected, which implies that $D\setminus\{p_0^1\}$ is not a dominating set.
	
	Suppose, for contradiction, that $G[D\setminus\{p_0^1\}]$ is not connected. Then there exists a vertex $d_j\in D$ whose only neighbor in $V_{P_0}\setminus\{p_0^{n_0}\}$ is $p_0^1$. Consequently, the vertices
	$
	v_j,d_j,p_0^1,p_0^2,\ldots,p_0^{n_0-1}
	$
	induce a chain in $G$ of order $n_0+1$, contradicting the fact that $P_0$ is a maximum induced chain of $G$.
	
	It follows that $D\setminus\{p_0^1\}$ is not a dominating set. Hence, there exists a vertex $u$ outside $D\setminus\{p_0^1\}$
	that is not adjacent to any vertex in this set. 
	Since $p_0^1$ and $p_0^{n_0}$ are adjacent to $p_0^2$ and $p_0^{n_0-1}$, respectively, $u$ cannot be either endpoint of $P_0$.
	 Moreover, $p_0^1$ is the only vertex in $D$ adjacent to $u$. Also, since the vertices $p_0^2,p_0^3,\ldots,p_0^{n_0-1}$ belong to $D$, $u$ is neither one of these vertices nor adjacent to any of them. Consequently, $u$ must be adjacent to $p_0^{n_0}$; otherwise, the vertices
	$
	u,p_0^1,p_0^2,\ldots,p_0^{n_0}
	$
	would induce a chain in $G$ of order $n_0+1$, contradicting the fact that $P_0$ is a maximum induced chain of $G$. Finally, note that $u$ is not adjacent to $d_1$, since $d_1\in D\setminus\{p_0^1\}$, and not to $v_1$, since $p_0^2,p_0^3,\ldots,p_0^{n_0},u,v_1
	$
	would induce a chain in $G$ of order $n_0+1$, contradicting the fact that $P_0$ is a maximum induced chain of $G$.
	
	Finally, let $q_0^1=p_0^2$, $q_0^2=p_0^3,\ldots,q_0^{n_0-1}=p_0^{n_0}$ and
	$
	q_0^{n_0}=u.
	$
	Then the vertices $q_0^1,q_0^2,\ldots,q_0^{n_0}$ induce a maximum chain $Q$ in $G$. Since \texttt{Build-Linear-Forest} could have started with $P_0=Q$ and subsequently constructed a second chain using the same vertices $v_1$ and $d_1$,
	and since the largest index $i$ for which $d_1$ is adjacent to a vertex $q_0^i$ of $Q$ is strictly smaller than the corresponding index for $P_0$, we obtain a contradiction to the choice of $P_0$ and $d_1$.
\end{proof}

\begin{corollary}\label{cor:L-1}
	Let $G$ be a connected graph of order at least $2$. If $\gamma_c(G)=L(G)-1$, then every maximum induced chain $P$ of $G$ has $|V_P|=L(G)$ vertices, and deleting either endpoint of $P$ results in a minimum connected dominating set of $G$.
\end{corollary}
\begin{proof}
	Let $P$ be a maximum induced chain in $G$ with vertex set $\{p_1,p_2,\ldots,p_{|V_P|}\}$. By Theorem~\ref{thm:dominating}, $V_P$ is a connected dominating set of $G$.
	Consider an execution of \texttt{Build-Linear-Forest} in which the first chain $P_0$ constructed at Step 6 is $P$. Since $V_P$ is a dominating set, the algorithm terminates with $P$ as its output. The associated connected dominating set $D$ consists of all vertices of $P$ except one of its endpoints. 
	Since $P$ is a chain, we have $|V_P|\leq L(G)$. Consequently,
	$$
	L(G)-1=\gamma_c(G)\leq |D|=|V_P|-1\leq L(G)-1.
	$$
	Equality must hold throughout, so $|V_P|=L(G)$. Moreover, $D$ has $\gamma_c(G)$ vertices and is therefore a minimum connected dominating set. Hence, deleting either endpoint of $P$ yields a minimum connected dominating set $D$ of $G$.
	\end{proof}

	\begin{corollary}\label{cor:common-neighbor}
	Let $G$ be a connected non-complete graph satisfying $\gamma_c(G)=L(G)-1$, and let $P$ be a maximum induced chain of $G$. Then there exists a vertex $q\notin V_P$ that is adjacent to the two endpoints of $P$ and to no other vertex of $P$.
\end{corollary}

\begin{proof}
	Let $P$ be a maximum induced chain of $G$, with vertices $p_1,p_2,\ldots,p_{|V_P|}$ ordered so that $p_j$ is adjacent to $p_{j+1}$ for $j=1,\ldots,|V_P|-1$. Let $I=\{p_2,\ldots,p_{|V_P|-1}\}$ be the vertex set obtained from $V_P$ by removing the two endpoints of $P$. 
	By Corollary \ref{cor:L-1}, the condition $\gamma_c(G)=L(G)-1$ implies that $|V_P|=L(G)$ and that both \(I\cup\{p_1\}\) and  \(I\cup\{p_{|V_P|}\}\) are connected dominating sets. Since $G$ is not complete, $L(G)=|V_P|\geq 3$, and hence $I\neq \emptyset$. Since $|I|=L(G)-2<\gamma_c(G)$, $I$ is not a dominating set, which means that there is a vertex $q\notin I$ with no neighbor in $I$. Since both sets \(I\cup\{p_1\}\) and \(I\cup\{p_{|V_P|}\}\) are dominating, $q$ must be adjacent to both $p_1$ and $p_{|V_P|}$.
\end{proof}

\section{The Case  $L(G)= \gamma_c(G)+1$}\label{sec:5}
As shown in the previous section, $L(G)\geq \gamma_c(G)+1$. We now show that, when $L(G)$ attains this lower bound, Conjecture \ref{conj:main} holds, assuming that it holds for graphs of smaller order.
\begin{theorem}\label{thm:gc-L-1}
	Let $G$ be a connected graph of order $n$ and $L(G)=\gamma_c(G)+1$. If Conjecture \ref{conj:main} holds for all connected graphs of order at most $n-1$, then it holds for $G$.
\end{theorem}

\begin{proof}
As shown in Theorem \ref{thm:L3}, Conjecture \ref{conj:main} holds when $L(G)\leq 3$. We may therefore assume that $L(G)\ge4$, so $G$ is not complete. By Corollary~\ref{cor:L-1}, $G$ contains an induced chain $P$ of order $L(G)$, and $V_P$ is a connected dominating set of $G$, since $V_P$ minus an endpoint of $P$ already is. Let $p_1,p_2,\dots,p_{L(G)}$ denote the vertices of $P$, ordered so that $p_j$ is adjacent to $p_{j+1}$ for $j=1,\dots,L(G)-1$. By Corollary~\ref{cor:common-neighbor}, there is a vertex $q\notin V_P$ whose only neighbors on $P$ are $p_1$ and $p_{L(G)}$. Since these two vertices are not adjacent, $d_G(p_1,p_{L(G)})=2$.

	Let $u$ be an endpoint of $P$. Suppose, for a contradiction, that $G-u$ is disconnected. Then there must be a vertex outside $P$ whose only neighbor in $P$ is $u$. Adding this vertex to $P$ would yield a linear forest of order $L(G)+1$, contradicting the definition of $L(G)$. Hence, $G-u$ is connected. By the induction hypothesis, $L(G-u)(n-1)(n-2)\geq 4W(G-u)$ for every endpoint $u$ of $P$.

	Define $\tau_G(u)=L(G)(n-1)-2w_G(u)$ for each endpoint $u$ of $P$, and $\sigma(v)=d_G(p_1,v)+d_G(p_{L(G)},v)$ for every vertex $v$ of $G$. Then
	\begin{equation}
		\tau_G(p_1)+\tau_G(p_{L(G)})=2L(G)(n-1)-2\sum_{v\in V_G}\sigma(v)=2\Bigl(\sum_{v\in V_G}\bigl(L(G)-\sigma(v)\bigr)-L(G)\Bigr).
	\end{equation}
	We bound the terms of the last sum.
	\begin{itemize}
		\item $\sigma(p_1)=\sigma(p_{L(G)})=\sigma(q)=2$, so $L(G)-\sigma(v)=L(G)-2\geq 2$ for $v\in\{p_1,q,p_{L(G)}\}$.
		\item For $2\le j\le L(G)-1$, the chain $P$ gives $\sigma(p_j)\le(j-1)+(L(G)-j)=L(G)-1$, so $L(G)-\sigma(p_j)\ge1$.
		\item Let $K$ be the set of vertices $v$ of $G$ such that $\sigma(v)>L(G)$; thus $L(G)-\sigma(v)\ge 0$ for every $v\notin K$. By the first two items, $K$ is disjoint from $V_P\cup\{q\}$. Moreover, every neighbor $v$ of $p_1$ satisfies $\sigma(v)\le1+\bigl(1+d_G(p_1,p_{L(G)})\bigr)=4\le L(G)$, and similarly for $p_{L(G)}$. So no vertex of $K$ is adjacent to $p_1$ or to $p_{L(G)}$.
	\end{itemize}
	Let $v\notin V_P$, and let $j\le j'$ be the minimum and maximum indices of its neighbors on $P$, which exist since $V_P$ is dominating. Then $d_G(p_1,v)\le j$ and $d_G(p_{L(G)},v)\le L(G)-j'+1$, so $\sigma(v)\le L(G)+1-(j'-j)$. Consequently, every vertex $v\in K$ has a unique neighbor $p_{j_v}$ on $P$, and $\sigma(v)=L(G)+1$. By the third item, $2\le j_v\le L(G)-1$. Furthermore, distinct vertices $v,v'\in K$ satisfy $j_v\neq j_{v'}$, since otherwise $(V_P\setminus\{p_{j_v}\})\cup\{v,v'\}$ would induce a linear forest of order $L(G)+1$, consisting of the two segments of $P-p_{j_v}$ together with $v$ and $v'$, which form either an edge or two isolated vertices. Therefore the map $v\mapsto p_{j_v}$ is an injection from $K$ into $\{p_2,\dots,p_{L(G)-1}\}$, and
	\begin{align*}
	\sum_{v\in V_G}\bigl(L(G)-\sigma(v)\bigr)
	&\ge 3(L(G)-2)+\sum_{j=2}^{L(G)-1}\bigl(L(G)-\sigma(p_j)\bigr)-\abs{K}\\
	&\ge 3(L(G)-2)+(L(G)-2)-\abs{K}
	\ \ge\ 3(L(G)-2).
	\end{align*}
	
	By \eqref{eq:tau-sum}, $\tau_G(p_1)+\tau_G(p_{L(G)})\ge2\bigl(3(L(G)-2)-L(G)\bigr)=4(L(G)-3)>0$.
	
	Hence at least one of $\tau_G(p_1)$ and $\tau_G(p_{L(G)})$ is positive; let $u$ denote the corresponding vertex. Every induced linear forest of $G-u$ is an induced linear forest of $G$, so $L(G-u)\le L(G)$. Distances in $G-u$ are at least those in $G$, so $W(G-u)\ge W(G)-w_G(u)$. Therefore
	\begin{align*}
		\delta(G)-\delta(G-u)
		&=L(G)n(n-1)-4W(G)-L(G-u)(n-1)(n-2)+4W(G-u)\\[4pt]
				&\ge L(G)n(n-1)-4W(G)-L(G)(n-1)(n-2)+4W(G)-4w_G(u)\\[4pt]
		&=2L(G)(n-1)-4w_G(u)=2\tau_G(u)>0.
	\end{align*}

	Since $G-u$ is connected, $\delta(G-u)\ge0$, and therefore $\delta(G)>0$.
\end{proof}

\section{The Case $L(G)\geq \max\{\gamma_c(G)+2,\diam(G)+3\}$}\label{sec:6}

As established in Theorem~\ref{thm:bound}, every connected graph $G$ satisfies $\gamma_c(G)\leq L(G)-1$. Moreover, Theorem~\ref{thm:gc-L-1} shows that Conjecture~\ref{conj:main} holds whenever $\gamma_c(G)=L(G)-1$.
Moreover, $L(G)\geq \diam(G)+1$ and we have shown in Theorems \ref{thm:diam1} and \ref{thm:diam2} that if the conjecture holds for connected graphs of order at most $n-1$ and if $L(G)\leq \diam(G)+2$, then it also holds for $G$. Hence, it remains to establish the conjecture for connected graphs $G$ satisfying $L(G)\geq \max\{\gamma_c(G)+2, \diam(G)+3\}$. In this section, we show that the inequality $L(G)n(n-1)-4W(G)\geq 0$ is strict whenever these two conditions are satisfied.

\begin{theorem}\label{thm:casdifficile}
	If $G$ is a connected graph such that $L(G)\geq \max\{\gamma_c(G)+2, \diam(G)+3\}$, then $L(G)n(n-1)-4W(G)>0$.
\end{theorem}

Before proceeding with the proof, we make a few useful observations.
Note that graphs satisfying $L(G)\geq \gamma_c(G)+2$ have $\diam(G)\geq 2$, since otherwise $G$ is complete and $L(G)=2<3=\gamma_c(G)+2$; hence $n\geq L(G)\geq 5$. 

So, let $G$ be a graph satisfying $L(G)\geq \max\{\gamma_c(G)+2, \diam(G)+3\}$. Let $D$ be a minimum connected dominating set of $G$. Let $T$ be a spanning tree of $G[D]$ and let $S$ be a spanning tree of $G$ obtained from $T$ by linking each vertex outside $T$ to one of its neighbors in $T$. Thus $T$ is a tree of order $t\leq L(G)-2$, while $S$ is of order $n$. By Property \ref{pro:Wiener}, we have $W(S)=\sum_{e\in E_S}\pi(S,e)$. The contribution $\pi(S,e)$ of each edge $e$ of $S$ to $W(S)$, is as follows:
\begin{itemize}
	\item For each of the $n-t$ edges $e$ added to $T$ to obtain $S$, we have $\pi(S,e)=n-1$, since deleting $e$ yields two connected components, one with a single vertex and the other with $n-1$ vertices.
	\item For each of the $t-1$ edges $e$ of $T$, we have $\pi(S,e)=\frac{n^2-a(S,e)^2}{4}$. Indeed, let $x$ be the order of the largest connected component of $S-e$. The other component then has $n-x$ vertices, so $a(S,e)=2x-n$ and
	$$
	\pi(S,e)=x(n-x)=\frac{n^2}{4}-\Big(x-\frac{n}{2}\Big)^2=\frac{n^2}{4}-\Big(\frac{2x-n}{2}\Big)^2=\frac{n^2-a(S,e)^2}{4}.
	$$
\end{itemize}
Hence, 
\begin{align}\label{4WS}4W(S)=4(n-t)(n-1)+(t-1)n^2-\sum_{e\in E_T}a(S,e)^2.
\end{align}
\begin{lemma}\label{case0}
	If $L(G)>\gamma_c(G)+2$, then $L(G)n(n-1)-4W(G)>0$.
	\end{lemma}
\begin{proof}
	In this case, we have $t=\gamma_c(G)\leq L(G)-3$. Since $(3t+1)n>4t$, using (\ref{4WS}), we obtain
	\begin{align*}
		&L(G)n(n-1)-4W(G)\\[0pt]
		\geq&L(G)n(n-1)-4W(S)\\[-2pt]
		=&L(G)n(n-1)-4(n-t)(n-1)-(t-1)n^2+\sum_{e\in E_T}a(S,e)^2\\[-3pt]
		=&(L(G)-t-3)n(n-1)+(3t+1)n-4t+\sum_{e\in E_T}a(S,e)^2>0.
	\end{align*}
\end{proof}
Henceforth, we assume that $\gamma_c(G)=t=L(G)-2$. Consequently, again using (\ref{4WS}), we have
\begin{align*}		
	&L(G)n(n-1)-4W(G)\\
	=&L(G)n(n-1)-4W(S)+4\Big(W(S)-W(G)\Big)\\
	=&(t+2)n(n-1)-4(n-t)(n-1)-(t-1)n^2+\sum_{e\in E_T}a(S,e)^2+4\Big(W(S)-W(G)\Big)\\[-4pt]
	=&{-}\Big(n^2{-}(3t{+}2)n{+}4t\Big){+}\sum_{e\in E_T}a(S,e)^2{+}4\Big(W(S)-W(G)\Big).
\end{align*}

We deduce that proving $L(G)n(n-1)-4W(G)> 0$ is equivalent to proving the following inequality:
\begin{align}
		\label{eq:SG}
		\sum_{e\in E_T}a(S,e)^2+4\Big(W(S)-W(G)\Big)> n^2-(3t+2)n+4t.
\end{align}

Since $W(S)\geq W(G)$, it suffices to prove 
\begin{align}\label{eq:SG2}
\sum_{e\in E_T}a(S,e)^2> n^2-(3t+2)n+4t.	
\end{align}

Moreover, for $h\ge \diam(G)$, we have 
			$$W(S)-W(G)=(W(S)-\underline{W}(S,h))+(\underline{W}(S,h)-W(G))$$

and it follows from Property \ref{prop:W} that both terms on the right-hand side of the equality are non-negative. Therefore, it suffices to prove
\begin{align}
	\label{eq:SG3}
	\sum_{e\in E_T}a(S,e)^2+4\Big(W(S)-\underline{W}(S,h)\Big)> n^2-(3t+2)n+4t.
\end{align}

The proof of Theorem \ref{thm:casdifficile} is divided into seven cases, depending on the number of leaves in $T$ and their configuration. In two cases (Lemmas~\ref{case1} and \ref{case5}), we establish inequality (\ref{eq:SG2}). In three cases (Lemmas \ref{case2}, \ref{case3}, and \ref{case4}), we establish inequality (\ref{eq:SG3}) with $h=t-1=L(G)-3\geq \diam(G)$. In another case (Lemma \ref{case61}), under the assumption that $L(G)\ge\diam(G)+4$, we establish inequality (\ref{eq:SG3}) with $h=t-2=L(G)-4\geq \diam(G)$. Finally,  in the remaining case (Lemma \ref{case6}), we establish  inequality (\ref{eq:SG}).

Before proceeding with the proofs, we make an observation that will be useful throughout. If $T$ contains $\lambda\geq 3$ edges $e_1,e_2,\ldots,e_{\lambda}$ such that each chain in $T$ contains at most 2 of these edges, then 
\begin{align}\label{eq:a}\sum_{i=1}^{\lambda}a(S,e_i)^2\geq n^2\frac{(\lambda-2)^2}{\lambda}.
\end{align}
Indeed, let $x_i$ denote the number of vertices in the connected component of $S-e_i$ containing no edge $e_j$ with $j\neq i$. By the Cauchy-Schwarz inequality, we have
$$
\sum_{i=1}^{\lambda}a(S,e_i)^2=\sum_{i=1}^{\lambda}|x_i-(n-x_i)|^2=\sum_{i=1}^{\lambda}(n-2x_i)^2\geq \frac{(\lambda n-2\sum_{i=1}^{\lambda}x_i)^2}{\lambda}\geq
\frac{(\lambda n-2n)^2}{\lambda}.$$

We can now proceed with the proof of Theorem \ref{thm:casdifficile}. We first consider in Lemmas \ref{case1} and \ref{case2} the cases in which the number of leaves of $T$ is not equal to 3.
	
\begin{lemma}\label{case1} If $T$ has at least four leaves then $$\sum_{e\in E_T}a(S,e)^2> n^2-(3t+2)n+4t.$$
\end{lemma}
\begin{proof}Let $e_1,e_2,e_3,e_4$ be four edges of $T$ each having an endpoint of degree 1. Applying (\ref{eq:a}) with $\lambda=4$, and noting that $(3t+2)n>4t$, we obtain  \begin{align*}\sum_{i=1}^4 a(S,e_i)^2\geq n^2>n^2-(3t+2)n+4t.&\qedhere\end{align*}
	\end{proof}
	
	\begin{lemma}\label{case2} If $T$ is a chain then $$\sum_{e\in E_T}a(S,e)^2+4\Big(W(S)-\underline{W}(S,t-1)\Big)> n^2-(3t+2)n+4t.$$
		\end{lemma}
		\begin{proof}
		Let $e_1$ and $e_2$ be the two edges of $T$ incident to its two leaves $f_1$ and $f_2$, and let $x_i$ denote the number of vertices in the connected component of $S-e_i$ that contains $f_i$. Hence, we have 
			\begin{align*}
			\sum_{e\in E_T}a(S,e)^2&\geq a(S,e_1)^2+a(S,e_2)^2=(2x_1-n)^2+(2x_2-n)^2\\&=4(x_1+x_2)^2-8x_1x_2-4n(x_1+x_2)+2n^2
\end{align*}
			Also, since 
			\begin{itemize}
				\item Each of the $x_1-1$ vertices of $S$ outside $T$ adjacent to  $f_1$ is at distance $t+1$ from each of the $x_2-1$ vertices of $S$ outside $T$ connected to $f_2$;
				\item $f_1$ is at distance at least $t$ from each of the $x_2-1$ vertices of $S$ outside $T$ connected to $f_2$, and $f_2$ is at distance $t$ from each of the $x_1-1$ vertices of $S$ outside $T$ connected to $f_1$,
			\end{itemize}
			we have
			$$
			W(S)-\underline{W}(S,t-1)\geq 2(x_1-1)(x_2-1)+(x_1-1)+(x_2-1)=2x_1x_2-x_1-x_2.
			$$
			Hence,  
			\begin{align*}	
				\sum_{e\in E_T}a(S,e)^2+4\Big(W(S)-\underline{W}(S,t-1)\Big)&\geq  4(x_1+x_2)^2-4(n+1)(x_1+x_2)+2n^2.
			\end{align*}
			This last function on the right-hand side is a quadratic in $y=x_1+x_2$ that attains its minimum at $y=\frac{n+1}{2}$. Therefore, since $3tn-4t-1>0$, we have
			$$\sum_{e\in E_T}a(S,e)^2+4\Big(W(S)-\underline{W}(S,t-1)\Big)\geq 
			(n+1)^2-2(n+1)^2+2n^2=n^2-2n-1>n^2-(3t+2)n+4t.
			$$
		\end{proof}
		
		We next consider the case in which $T$ has three leaves, precisely one of which is adjacent to the vertex of degree 3.
		
	\begin{lemma}\label{case3} If $T$ has exactly three leaves, one at distance 1 from the vertex of degree 3, and the two others at distance at least 2, then $$\sum_{e\in E_T}a(S,e)^2+4\Big(W(S)-\underline{W}(S,t-1)\Big)> n^2-(3t+2)n+4t.$$
\end{lemma}
\begin{proof}
	Let $f_1$ and $f_2$ be the two leaves at distance at least 2 from the vertex of degree 3, and let $f_3$ be the leaf adjacent to it. Let $e_1$, $e_2$ be  the two edges of $T$ incident to $f_1$ and $f_2$, and let $x_i$ denote the number of vertices in the connected component of $S-e_i$ that contains $f_i$. Also, let $e_3, e_4, e_5$ be the three edges incident to the vertex of degree 3 in $T$. Note that the five edges $e_i$, $i=1,2,\ldots 5$ are distinct. Applying (\ref{eq:a}) with the $\lambda=3$ edges $e_3,e_4,e_5$ gives
\begin{align*}
	\sum_{e\in E_T}a(S,e)^2&\geq \sum_{i=1}^5a(S,e_i)^2 \geq (2x_1-n)^2+(2x_2-n)^2+\frac{n^2}{3}\\&= 4(x_1+x_2)^2-8x_1x_2-4n(x_1+x_2)+\frac{7n^2}{3}.
\end{align*}	Each of the $x_1-1$ vertices of $S$ outside $T$ connected to $f_1$ is at distance at least $t$ from each of the $x_2{-}1$ vertices of $S$ outside $T$ connected to $f_2$, which implies
	$$
	W(S)-\underline{W}(S,t-1)\geq (x_1-1)(x_2-1)=x_1x_2-(x_1+x_2)+1.
	$$
	Using $4x_1x_2\leq (x_1+x_2)^2$, we get
	\begin{align*}	
		\sum_{e\in E_T}a(S,e)^2+4\Big(W(S)-\underline{W}(S,t-1)\Big)&\geq 4(x_1+x_2)^2-4x_1x_2-4(n+1)(x_1+x_2)+\frac{7n^2}{3}+4\\[-8pt]
		&\geq 3(x_1+x_2)^2-4(n+1)(x_1+x_2)+\frac{7n^2}{3}+4
	\end{align*}
	This last function on the right-hand side is a quadratic in $y=x_1+x_2$ that attains its minimum at $y=\frac{2(n+1)}{3}$. Therefore, since $n(9t-2)-12t+8>0$, we have
	\begin{align*}\sum_{e\in E_T}a(S,e)^2+4\Big(W(S)-\underline{W}(S,t-1)\Big)&\geq 
		\frac{4}{3}(n+1)^2-\frac{8}{3}(n+1)^2+\frac{7}{3}n^2+4\\[-7pt]&=n^2-\frac{8}{3}n+\frac{8}{3}>n^2-(3t+2)n+4t.
	\end{align*}
\end{proof}
		
		We next consider the case in which $T$ has exactly three leaves, at least two of which are adjacent to the vertex of degree 3.
	\begin{lemma}\label{case4} If $T$ has exactly three leaves, at least two of which are adjacent to the vertex of degree 3, then $$\sum_{e\in E_T}a(S,e)^2+4\Big(W(S)-\underline{W}(S,t-1)\Big)> n^2-(3t+2)n+4t.$$
\end{lemma}
\begin{proof}
	Let $f_1$ and $f_2$ be the two leaves at distance 1 from the vertex of degree 3, and let $f_3$ be the third leaf. For $i=1,2,3$, let $e_i$ be  the edge of $T$ incident to $f_i$, and let $x_i$ denote the number of vertices in the connected component of $S-e_i$ that contains $f_i$. Each of the $x_1-1$ vertices of $S$ outside $T$ connected to $f_1$ and each of the $x_2-1$ vertices of $S$ outside $T$ connected to $f_2$ is at distance $t$ from each of the $x_3-1$ vertices of $S$ outside $T$ connected to $f_3$, which implies
\begin{align*}
	\sum_{e\in E_T}a(S,e)^2+4\Big(W(S)-\underline{W}(S,t-1)\Big)\geq \sum_{i=1}^2(2x_i-n)^2+(2x_3-n)^2+4(x_1+x_2-2)(x_3-1).
\end{align*}	
The last two terms form a quadratic in $x_3$ whose minimum, attained at $x_3=\frac{n-x_1-x_2+2}{2}$, is $(x_1+x_2-2)(2n-x_1-x_2-2)$. By the Cauchy-Schwarz inequality, we have
$(2x_1-n)^2+(2x_2-n)^2\geq 2(n-x_1-x_2)^2$. Hence,
writing $y=x_1+x_2$, we obtain
\begin{align*}
\sum_{e\in E_T}a(S,e)^2+4\Big(W(S)-\underline{W}(S,t-1)\Big)&\geq 2(n-y)^2+(y-2)(2n-y-2)\\[-6pt]
&=(n-y)^2+n^2-4n+4\\&> n^2-4n+4>n^2-(3t+2)n+4t.\qedhere
\end{align*}	
	\end{proof}
We now consider the case where $T$ has exactly three leaves, each at distance at least 3 from the vertex of degree 3.		\begin{lemma}\label{case5} If $T$ has exactly three leaves, all of them at distance at least 3 from the vertex of degree 3, then $$\sum_{e\in E_T}a(S,e)^2> n^2-(3t+2)n+4t.$$
	\end{lemma}
	\begin{proof}
		Let $\ell_1$, $\ell_2$, and $\ell_3$ denote the lengths of the three chains in $T$ joining the vertex of degree 3 to the three leaves. By assumption, we have $\frac{1}{\ell_1}+\frac{1}{\ell_2}+\frac{1}{\ell_3}\leq 1$. For $i=1,2,3$, let $e_i$ be the edge of the $i$-th chain incident to the vertex of degree $3$, and $x_i$ the number of vertices of the connected component of $S-e_i$ not containing that vertex; thus $x_1+x_2+x_3\le n-1$. For every edge $e$ of the $i$-th chain, the component of $S-e$ not containing the vertex of degree $3$ has at most $x_i$ vertices, so $a(S,e)^2\ge(n-2x_i)^2$ whenever $x_i\le n/2$. If $x_1,x_2,x_3\le n/2$, the Cauchy--Schwarz inequality gives
		\begin{align*}
			\sum_{e\in E(T)}a(S,e)^2
			&\ge\sum_{i=1}^3\ell_i(n-2x_i)^2
			\ge\frac{\bigl(3n-2(x_1+x_2+x_3)\bigr)^2}{1/\ell_1+1/\ell_2+1/\ell_3}
			\ge(n+2)^2.
		\end{align*}
		Otherwise, one of them, say $x_1$, exceeds $n/2$; then $x_2+x_3<n/2$, and the same argument gives
				\begin{align*}
			\sum_{e\in E(T)}a(S,e)^2
			&\ge\sum_{i=2}^3\ell_i(n-2x_i)^2
			\ge\frac{\bigl(2n-2(x_2+x_3)\bigr)^2}{1/\ell_2+1/\ell_3}
			> n^2.
		\end{align*}
Hence, in both cases, $\sum_{e\in E_T}a(S,e)^2> n^2>n^2-(3t+2)n+4t.$
		\end{proof}
		
It remains to consider the situation in which $T$ has exactly three leaves, one at distance 2 from the vertex of degree 3 and the two others at distance at least 2 from it. We first consider the case where $L(G)\geq \diam(G)+4$.
	\begin{lemma}\label{case61}  If $L(G)\geq \diam(G)+4$ and $T$ has exactly three leaves, one at distance 2 from the vertex of degree 3, and the two others at distance at least 2, then $$\sum_{e\in E_T}a(S,e)^2+4\Big(W(S)-\underline{W}(S,t-2)\Big)> n^2-(3t+2)n+4t.$$
\end{lemma}
\begin{proof}
	In this case, $T$ has two leaves $f_1$ and $f_2$ at distance $t-3$ from each other and a third leaf $f_3$ at distance 2 from the vertex of degree 3. Let $e_1$, $e_2$ be  the two edges of $T$ incident to $f_1$ and $f_2$, and let $x_i$ denote the number of vertices in the connected component of $S-e_i$ that contains $f_i$. Also, let $e_3, e_4, e_5$ be the three edges incident to the vertex of degree 3 in $T$. Note that the five edges $e_i$, $i=1,2,\ldots 5$ are distinct. As in Lemma \ref{case3}, applying (\ref{eq:a}) to $e_3,e_4,e_5$ with $\lambda=3$ gives
	$$
	\sum_{e\in E(T)}a(S,e)^2\geq 
	(2x_1-n)^2+(2x_2-n)^2+\frac{n^2}{3}.
	$$
	Each of the $x_1-1$ vertices of $S$ outside $T$ connected to $f_1$ is at distance at least $t-1$ from each of the $x_2{-}1$ vertices of $S$ outside $T$ connected to $f_2$, which implies
	$$
	W(S)-\underline{W}(S,t-2)\geq (x_1-1)(x_2-1).
	$$
	Using the same calculations as in Lemma~\ref{case3}, we obtain
		\begin{align*}\sum_{e\in E_T}a(S,e)^2+4\Big(W(S)-\underline{W}(S,t-2)\Big)&\geq n^2-\frac{8}{3}n+\frac{8}{3}>n^2-(3t+2)n+4t.\qedhere
	\end{align*}
\end{proof}
		For the remainder of the proof, we may therefore assume that $L(G)=\diam(G)+3$. We will prove inequality (\ref{eq:SG}) by establishing lower bounds for both $\sum_{e\in E_T}a(S,e)^2$ and $W(S)-W(G)$. For bounding $W(S)-W(G)=\sum_{\{u,v\}\subseteq V_G}(d_S(u,v)-d_G(u,v))$, we focus on pairs of adjacent vertices $u,v$ in $G$ that lie outside $T$. Indeed, for each such pair, $d_G(u,v)=1$, whereas $d_S(u,v)\geq 2$. Consequently, each of these pairs contributes at least one unit to $W(S)-W(G)$.
		
\begin{figure}[H]
			\centering
			\includegraphics[width=0.74\linewidth]{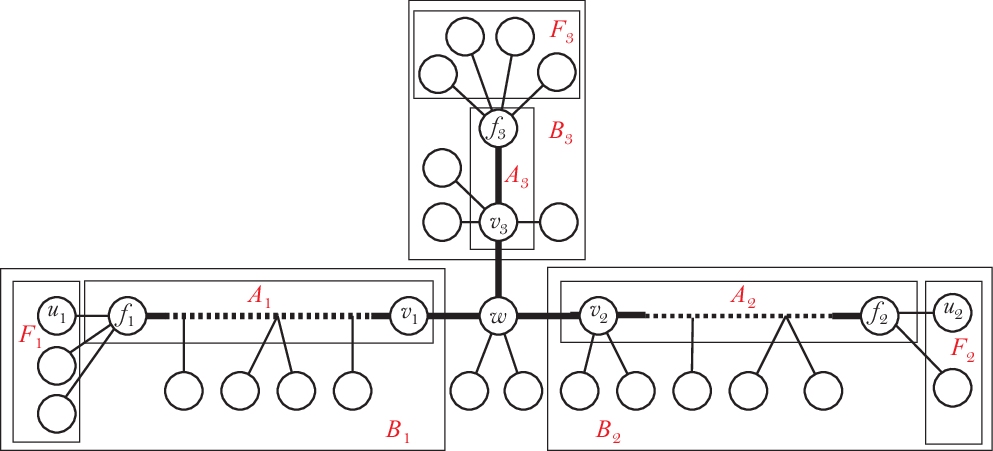}
			\caption{the tree $S$ and its subtree $T$ under the conditions stated in Lemma \ref{case6}.}
			\label{fig:ST}
		\end{figure}

		Since the structure of $T$ and its extension to $S$ are not arbitrary, we first introduce some notation, all of which are illustrated in Figure \ref{fig:ST}. In this figure, the edges of $T$ are drawn with thick lines, whereas the edges of $S$ connecting a vertex outside $T$ to a vertex in $T$ are drawn with thin lines. As mentioned, we assume that $T$ has exactly three leaves and denote by $w$ its vertex of degree 3. Its three neighbors in $T$ are denoted by $v_1$, $v_2$, and $v_3$. The three leaves of $T$ are denoted by $f_1$, $f_2$, and $f_3$, where $f_i$ belongs to the component containing $v_i$ after the three edges incident to $w$ are removed from $T$. We assume that $f_3$ is at distance 2 from $w$, meaning that $f_3$ and $v_3$ are adjacent, whereas $f_1$ and $f_2$ are at distance at least 2 from $w$. We denote by $A_i$ the set of vertices on the chain that connects $v_i$ to $f_i$ in $T$,  whereas $B_i$ is the vertex set of the connected component containing $f_i$ in the graph obtained from $S$ by removing the edge joining $w$ to $v_i$. Finally, for $i=1,2,3$, we denote by $F_i$ the set of vertices outside $T$ that are adjacent to $f_i$ in $S$.
		
			\begin{lemma}\label{case6} If $G$ is a connected graph of order $n$ such that $L(G)=\diam(G)+3=\gamma_c(G)+2$, and if $T$ has exactly three leaves, one at distance 2 from the vertex of degree 3, and the two others at distance at least 2, then  $$\sum_{e\in E_T}a(S,e)^2+4\Big(W(S)-W(G)\Big)> n^2-(3t+2)n+4t.$$
		\end{lemma}
		\begin{proof}		
		First, we note that we have some flexibility in the choice of $T$ and its extension to $S$. Indeed, we may choose any minimum connected dominating set $D$, and then take $T$ to be any spanning tree of $G[D]$. In what follows, whenever possible, we assume that $T$ is chosen so as to satisfy the conditions of Lemmas \ref{case1} through \ref{case5}. In other words, in the remaining case, we assume that no tree $T$ satisfying the conditions stated in the preceding lemmas can be obtained. When extending $T$ to $S$, a vertex $x$ outside $T$ may have several neighbors in $T$. To connect $x$ to $T$ in $S$, we give priority to neighbors other than $f_1$, $f_2$, $f_3$, and $v_3$. If none of these is available, we prefer $f_1$ or $f_2$ to $f_3$ or $v_3$. As a third priority, we prefer $v_3$ to $f_3$. Finally, if $f_3$ is the only neighbor of $x$ in $T$, then $x$ is connected to $f_3$ in $S$.
		
		Let $u_1$ and $u_2$ be two vertices such that $d_G(u_1,u_2)=\diam(G)$. Since $T$ has $t=L(G)-2=\diam(G)+1$ vertices, any two vertices of $T$ are at distance at most $t-3=\diam(G)-2$ in $T$. Hence, both vertices $u_1$ and $u_2$ lie outside $T$. Moreover, one of them, say $u_1$, is connected to $f_1$ in $S$, while the other, $u_2$, is connected to $f_2$. Consequently, the subgraph induced by $A_1\cup A_2\cup \{u_1,u_2,w\}$  is a chain of $G$, and $\diam(G)\geq 6$.
		
				If a vertex $x$ in $G$ is adjacent to both $f_1$ and $f_3$, then, for every $y_1\in F_1$ and every $y_3\in F_3$, we have $d_S(y_1,y_3)=|A_1|+4\geq 6$, whereas $d_G(y_1,y_3)\leq 4$, since $y_1,f_1,x,f_3,y_3$ forms a chain of length 4 in $G$. Thus, each such pair $y_1,y_3$ contributes at least two units to $W(S)-W(G)$. Let $e_1$ and $e_3$ be the edges of $T$ incident to $f_1$ and $f_3$, respectively. Their contribution to $\sum_{e\in E_T} a(S,e)^2$ is $(n{-}2(|F_1|{+}1))^2{+}(n{-}2(|F_3|{+}1))^2$. Since $|F_1|+|F_3|+2<n$ and $8-8n>4t-(3t+2)n$, we have, 
		\begin{align*}
			\sum_{e\in E_T} a(S,e)^2+4(W(S)-W(G))\geq&(n-2\big(|F_1|+1)\big)^2+(n-2\big(|F_3|+1)\big)^2+8|F_1||F_3|\\[-8pt]
			=&n^2+\big(n-2(|F_1|+|F_3|+2)\big)^2-8(|F_1|+|F_3|+2)+8\\
			>&n^2-8n+8>n^2-(3t+2)n+4t.
		\end{align*}
		The same argument applies with \(f_2\) in place of \(f_1\). We may therefore assume, in what follows, that no vertex is adjacent to \(f_3\) and to at least one of \(f_1\) and \(f_2\). In particular, since $u_1$ is adjacent to $f_1$ and $u_2$ to $f_2$, we deduce that $f_3$ is adjacent to neither $u_1$ nor $u_2$.
		
		Suppose that $B_3$ contains more than $\frac{n}{2}$ vertices. Then $|B_1|+|B_2|<\frac{n}{2}-1$. For every edge $e$  on the chain in $T$ from $w$ to $f_1$, we have $a(S,e)\geq n-2|B_1|$, since the connected component of $S-e$ containing $f_1$ has at most $|B_1|$ vertices. Similarly, $a(S,e)\geq n-2|B_2|$ for every edge $e$ on the chain from $w$ to $f_2$. Since  $|A_1|\geq 2$ and $|A_2|\ge 2$, the Cauchy–Schwarz inequality yields $$\sum_{e\in E_T}a(S,e)^2\ge 2(n-2|B_1|)^2+2(n-2|B_2|)^2\ge (2n-2(|B_1|+|B_2|))^2>(n+2)^2>n^2-(3t+2)n+4t.$$
		
		We may therefore assume $|B_3|\le \frac{n}{2}$. We now examine the possible neighbors of $v_3$ and $f_3$ in $G$. 
		\begin{itemize}
			\item If $v_3$ is adjacent to $f_1$ in $G$, replacing the edge joining $v_3$ to $w$ with the edge joining $v_3$ to $f_1$ yields a tree $T'$ that is a chain and therefore satisfies the conditions of Lemma~\ref{case2}, contradicting our assumption that $T$ was chosen so that no such tree could be obtained. The same argument applies if $v_3$ is adjacent to $f_2$. Hence, $N^1_G(v_3)\cap\{f_1,f_2\}=\emptyset$.
			\item If $f_3$ is adjacent to a vertex $x$ on $A_1\cup A_2\cup\{w\}$, replacing the edge joining $f_3$ to $v_3$ with the edge joining $f_3$ to $x$ yields a tree $T'$ with four leaves (if $x\neq f_1,f_2$), or a tree with exactly three leaves, one at distance $1$ from $w$ and the other two at distance at least $2$ from $w$. Hence $T'$ satisfies the conditions of  Lemma~\ref{case1} or~\ref{case3}, a contradiction. Therefore, $ N^1_G(f_3)\cap V_T=\{v_3\}$.
			\item Vertex $f_3$ must have at least one neighbor outside $T$ in $S$, since otherwise our priority rules for extending $T$ to $S$ would force $v_3$ to be the only neighbor of $f_3$ in $G$, making $V_T\setminus\{f_3\}$ a connected dominating set of size $\gamma_c(G)-1$, a contradiction. Hence, $|F_3|\geq 1$.
		\end{itemize} 
		
		Let $C=A_1\cup A_2\cup\{w\}$ be the set of vertices on the chain linking $f_1$ to $f_2$ in $T$. Our priorities for
		extending $T$ to $S$ ensure that no vertex in $B_3\setminus\{f_3,v_3\}$ is adjacent in $G$ to a vertex in $C$. Since $v_3$ is the only vertex in $T$ adjacent to $f_3$ in $G$, we deduce that $N_G^1(B_3\setminus\{v_3\})\cap C=\emptyset$.
		 Moreover, at least one of the sets $N_G^1(u_1)\cap (B_3\setminus\{v_3\})$ and $N_G^1(u_2)\cap (B_3\setminus\{v_3\})$ is empty. Indeed, $d_G(r,r')\leq 3$ for every pair of vertices $r,r'\in B_3\setminus\{v_3\}$, and if both sets were nonempty, we would have $d_G(u_1,u_2)\leq 5<\diam(G)$. Therefore,
		\begin{itemize}
			\item if both sets are empty, adding a linear forest of $G[B_3\setminus\{v_3\}]$ to the chain $G[C\cup\{u_1,u_2\}]$ yields a linear forest of $G$ with $t+L(G[B_3\setminus\{v_3\}])$ vertices. Since $L(G)=t+2$, this implies that $L(G[B_3\setminus\{v_3\}])\leq 2$. Moreover, since $B_3\setminus\{v_3\}$ contains $f_3$ and at least one neighbor in $F_3$,  $G[B_3\setminus\{v_3\}]$ is a complete graph. Hence, $G[B_3\setminus \{f_3,v_3\}]$ has $\binom{|B_3|-2}{2}$ edges.
			\item if one of the sets, say $N_G^1(u_1)\cap (B_3\setminus\{v_3\})$, is non empty, then adding a linear forest of $G[B_3\setminus\{v_3\}]$ to $G[C\cup\{u_2\}]$  yields a linear forest of $G$ with $t-1+L(G[B_3\setminus\{v_3\}])$ vertices. Since $L(G)=t+2$, this implies that $L(G[B_3\setminus\{v_3\}])\leq 3$. It follows from Property \ref{prop:nonedge} that $G[B_3\setminus \{f_3,v_3\}]$ has at least $\binom{|B_3|-2}{2}-(|B_3|-1)=\frac{|B_3|^2-7|B_3|+8}{2}$ edges.
		\end{itemize}
Let us now consider the subgraph $G[(C\cup\{u_2,f_3\})\setminus\{f_1\}]$ of $G$. It consists of the chain joining the neighbor of $f_1$ in $T$ to $u_2$, together with the isolated vertex $f_3$. It is therefore a linear forest on $L(G)-3$ vertices. The vertices in $F_1$ are not adjacent to $f_3$, since we have established that no vertex can be adjacent to both $f_1$ and $f_3$. Moreover, no vertex $x\in F_1$ is adjacent to $u_2$, since the chain induced by $u_1,f_1,x,u_2$ would imply that $\diam(G)=d_G(u_1,u_2)\leq 3$, contradicting the fact that $\diam(G)\geq 6$. Finally, by the priority rules defined for extending $T$ to $S$, no vertex in $F_1$ is adjacent to any vertex in $C\setminus \{f_1\}$. We conclude that any linear forest in $G[F_1]$ can be added to the linear forest $G[(C\cup\{u_2,f_3\})\setminus\{f_1\}]$ of $L(G)-3$ vertices. It follows that $L(G[F_1])\leq 3$ and we deduce from Property \ref{prop:nonedge} that $G[F_1]$ contains at least $\binom{|F_1|}{2}-(|F_1|)=\frac{|F_1|(|F_1|-3)}{2}$ edges. Similarly, $G[F_2]$ contains at least $\frac{|F_2|(|F_2|-3)}{2}$ edges.

For $i=1,2,3$, let $e_i^1$ be the edge incident to $f_i$ in $T$ and let $e_i^2$ be the edge incident to $w$ on the chain linking $w$ to $f_i$. Also, let $m_1,m_2$ and $m_3$ denote the number of edges in the subgraphs $G[F_1]$, $G[F_2]$ and $G[B_3\setminus \{f_3,v_3\}]$, respectively. Each edge in these three subgraphs contributes at least one unit to $W(S)-W(G)$, since it joins two vertices $x$ and $y$ such that $d_G(x,y)=1$ while $d_S(x,y)\geq 2$. Hence, 
$$\sum_{e\in E_T}a(S,e)^2+4(W(S)-W(G))\ge \sum_{i=1}^3\Big(a(S,e_i^1)^2+a(S,e_i^2)^2+4m_i\Big).$$
We now bound the right-hand side by splitting it into four terms.
\begin{itemize}
	\item We have shown that $m_1\geq \frac{|F_1|(|F_1|-3)}{2}$. Hence,
	\begin{align*}
		a(S,e_1^1)^2+4m_1\ge&\Big(n-2(|F_1|+1)\Big)^2+2|F_1|(|F_1|-3)\\
		=&6|F_1|^2-|F_1|(4n-2)+(n^2-4n+4).
	\end{align*}
	The minimum of this quadratic is attained at $|F_1|=\frac{2n-1}{6}$, which gives
	\begin{align}\label{borne1} a(S,e_1^1)^2+4m_1\ge \frac{2n^2-20n+23}{6}.
	\end{align}
\item	Similarly, 
\begin{align}\label{borne2} a(S,e_2^1)^2+4m_2\ge \frac{2n^2-20n+23}{6}.
\end{align}	
\item By the Cauchy–Schwarz inequality, we have $$a(S,e_1^2)^2+a(S,e_2^2)^2=(n-2|B_1|)^2+(n-2|B_2|)^2\ge 2(n-(|B_1|+|B_2|))^2.$$
Since $|B_1|+|B_2|+|B_3|\le n-1$, we have
\begin{align}\label{borne3}
	a(S,e_1^2)^2+a(S,e_2^2)^2\geq 2(|B_3|+1)^2
\end{align}
\item We have shown that $m_3\geq \frac{|B_3|^2-7|B_3|+8}{2}$. Hence,
\begin{align}\label{borne4}
	a(S,e_3^1)^2+a(S,e_3^2)^2+4m_3\geq&\Big(n-2(|B_3|-1)\Big)^2+\Big(n-2|B_3|\Big)^2+2(|B_3|^2-7|B_3|+8)\nonumber\\
	=&10|B_3|^2-(8n+22)|B_3|+(2n^2+4n+20).
\end{align}
\end{itemize}
Summing the four bounds given in (\ref{borne1}) --(\ref{borne4}), we obtain
$$\sum_{e\in E_T}a(S,e)^2+4(W(S)-W(G))\ge \frac{1}{3}\Big(4\big(3|B_3|-n\big)^2+4n^2-8n-54|B_3|+89\Big).$$
We have established that we may assume $|B_3|\leq n/2$, which gives 
$$\sum_{e\in E_T}a(S,e)^2+4(W(S)-W(G))\ge \frac{1}{3}\Big(4n^2-35n+89\Big).$$
The difference between the right-hand side of the last inequality and the desired lower bound  $n^2-(3t+2)n+4t$ is $\frac{1}{3}\big(n^2+(9t-29)n+(89-12t)\big)$, which is strictly positive whenever $t\geq 7$ and $n\geq t+3$. Hence, 
\begin{align*}\sum_{e\in E_T}a(S,e)^2+4(W(S)-W(G))>n^2-(3t+2)n+4t.&\qedhere\end{align*}
		\end{proof}
We now provide an overview of the proof of Theorem \ref{thm:casdifficile}.\\

\emph{Proof} of Theorem \ref{thm:casdifficile}
	
	Let $G$ be a graph of order $n$ satisfying $L(G)\geq \max\{\gamma_c(G)+2,\diam(G)+3\}$. We have to prove that  $L(G)n(n-1)-4W(G)>0$. Lemma \ref{case0} shows that the inequality is satisfied in the case where $L(G)>\gamma_c(G)+2$. For $L(G)=\gamma_c(G)+2$, we considered a dominating tree $T$ on $\gamma_c(G)$ vertices and established that the inequality holds in each of the following cases:
	\begin{itemize}
		\item  $T$ has at least four leaves: the inequality is  established in Lemma \ref{case1}.
		\item $T$ has two leaves: the inequality is  established in Lemma \ref{case2}.
		\item $T$ has three leaves, one at distance 1 from the vertex of degree 3, and the two others at distance at least 2: the inequality is  established in Lemma \ref{case3}.
		\item $T$ has three leaves, at least two at distance 1 from the vertex of degree 3: the inequality is  established in Lemma \ref{case4}.
		\item $T$ has three leaves, all of them at distance at least 3 from the vertex of degree 3: the inequality is  established in Lemma \ref{case5}.		
		\item $T$ has three leaves, one at distance 2 from the vertex of degree 3, and the two others at distance at least 2 from it: 
		\begin{itemize}
			\item if $L(G)\geq \diam(G)+4$, the inequality is 
			established in Lemma \ref{case61}	\item if $L(G)= \diam(G)+3$, the inequality is 
			established in Lemma \ref{case6}.\end{itemize}
		\end{itemize}	
		Since these cases cover all possible configurations, the proof of the theorem is complete.

\section{Conclusion}
We have established that Conjecture \ref{conj:main} stating that $\mu(G) \leq \frac{L(G)}{2}$ holds for every connected graph $G$. To prove this result, we first showed in Theorem \ref{thm:L3} that the conjecture holds whenever $L(G)\leq 3$. Then, using the fact that $L(G) \geq \diam(G)+1$, we proved in Theorems \ref{thm:diam1} and \ref{thm:diam2} that the conjecture holds when $L(G)=\diam(G)+1$ and $L(G)=\diam(G)+2$, respectively. We subsequently established that $L(G)\geq \gamma_c(G)+1$ and proved in Theorem \ref{thm:gc-L-1} that the conjecture holds when $L(G)=\gamma_c(G)+1$. Finally, Theorem \ref{thm:casdifficile} addresses the only remaining case, namely, $L(G)\geq\max\{\gamma_c(G)+2,\diam(G)+3\}$.

\end{document}